# Structure Calculation and Reconstruction of Discrete State Dynamics from Residual Dipolar Couplings using REDCRAFT

Casey A. Cole[1], Rishi Mukhapadhyay[1], Hanin Omar[1], Mirko Hennig[2], Homayoun Valafar[1,*]

[1]Department of Computer Science & Engineering, University of South Carolina, Columbia, SC 29208, USA

[2]Nutrition Research Institute, University of North Carolina at Chapel Hill, Kannapolis, NC 27514, USA

[*] Corresponding Author Email: homayoun@cec.sc.edu Phone: 1 803 777 2880 Fax: 1 803 777 3767

# Table of Contents

# 1 Abstract

Residual Dipolar Couplings (RDCs) acquired by Nuclear Magnetic Resonance (NMR) spectroscopy can be an indispensable source of information in investigation of molecular structures and dynamics. Here we present a complete strategy for structure calculation and reconstruction of discrete state dynamics from RDC data. Our method utilizes the previously presented REDCRAFT software package and its dynamic-profile analysis to complete the task of fragmented structure determination and identification of the onset of dynamics from RDC data. Fragmented structure determination was used to demonstrate successful structure calculation of static and dynamic domains for several models of dynamics. We provide a mechanism of producing an ensemble of conformations for the dynamical regions that describe any observed order tensor discrepancies between the static and dynamic domains within a protein. In addition, the presented method is capable of approximating relative occupancy of each conformational state.

The developed methodology has been evaluated on simulated RDC data with ±1Hz of error from an 83 residue α protein (PDBID 1A1Z), and a 213 residue α/β protein DGCR8 (PDBID 2YT4). Using 1A1Z, various models of arc and complex two and three discrete-state dynamics were simulated. MD simulation was used to generate a 2-state dynamics for DGCR8. In both instances our method reproduced structure of the protein including the conformational ensemble to within less than 2Å. Based on our investigations, arc motions with more than 30° of rotation are recognized as internal dynamics and are reconstructed with sufficient accuracy. Furthermore, states with relative occupancies above 20% are consistently recognized and reconstructed successfully. Arc motions with magnitude of 15° or relative occupancy of less than 10% are consistently unrecognizable as dynamical regions.

## 2 Introduction

Structural elucidation, including the study of molecular complexes and characterization of internal dynamics of macromolecules, is often the requisite step in the molecular characterization of diseases. Therefore, the development of rapid, cost-effective, automated, and unbiased methods of structure elucidation and study of internal dynamics of biological macromolecules becomes critical in the comprehension and treatment of many diseases. Our understanding of diseases at the molecular level is often impeded when encountering exceptionally challenging, yet important proteins such as membrane bound or associated proteins. These challenging classes of proteins are refractory for structure determination by traditional approaches. In the case of X-ray crystallography, this can be attributed to difficulty in obtaining a diffraction-quality crystal, whereas NMR spectroscopy is often hampered by an insufficient number of distance constraints.

Dynamical proteins are even more difficult to study which is directly related to their internal motion. Functionally relevant events often take place on time-scales (μsec-msec) that are inherently difficult to characterize. However, mounting evidence suggests that internal dynamics of proteins play an important role in their biological function. The breathing motion of myoglobin[1–4] can be cited as a historical instance of this property while other proteins such as hTS[5,64], DHFR[7–9], carbohydrate recognition domain of Galectin-3[10] can be cited as examples of dynamical proteins that are targets of pharmaceutical developments. Many RNA-binding proteins such as DGCR8, an integral component of the MicroRNA processing machinery, undergo conformational changes to enable the biological function of that protein. Therefore development of methods leading to elucidation of structure and dynamics of proteins is of paramount importance.

At present, comprehensive characterization of inter-molecular dynamics of a protein (or their complexes) has been limited to studies conducted by NMR spectroscopy. These studies are typically performed in successive steps beginning with the protein's structure determination under the assumption of rigidity and followed by characterization of its motion. Although structure determination protocols based on the assumption of molecular rigidity will yield a single structure, the degree of similarity between the static model of a protein structure and the many conformations sampled employing the dynamic model remains poorly understood. Our recent work[11,12] highlighted the possibility of obtaining an erroneous structure for a protein undergoing dynamics when disregarding flexibility during the structure determination step. Similarly, mapping of dynamics onto a false, static structure may lead to a compromised model of motion. This can be attributed to the fact that it is conceptually difficult to separate structure from dynamics, because the two are intimately related. In fact, any attempt in structure elucidation that disregards the dynamics of the protein (or vice versa), may produce faulty results. In addition, the currently existing methods of investigating internal dynamics of macromolecules by NMR spectroscopy are ineffective in sampling the biologically important time-scales ranging from picoseconds to millisecond[13–15].

A more practical and rigorous approach is to simultaneously characterize a protein's structure and its dynamics. Simultaneous characterization of structure and dynamics, and investigation of dynamics on biologically relevant time-scale are feasible through the use of Residual Dipolar Couplings[12,16]. Residual Dipolar Couplings (RDCs) are a class of experimental data acquired by NMR spectroscopy that may lead to the development of rapid, cost effective, automated and comprehensive approaches to characterization of structure and dynamics for both routine and challenging macromolecules[11,17–22]. RDCs already play a crucial role in providing reliable approaches for structure elucidation and possess the potential for integrating structure determination by NMR spectroscopy[11,16,23–25], X-ray crystallography[26–28], and computational modeling[26,27,29]. Importantly, RDC data can also be used as a probe of internal dynamics on biologically relevant time-scales[30,31].

Historically, the use of RDCs has been limited by two factors: data acquisition, and data analysis. The introduction of a wider variety of alignment media, combined with advances in instrumentation and data acquisition continue to address the experimental limitations in obtaining RDCs. The major bottleneck in utilization of RDC data in recent years has been attributed to a lack of RDC analysis tools capable of extracting the pertinent information embedded within this complex source of data. Nearly all legacy NMR data analysis software packages (i.e. Xplor-NIH[32], CNS[33]) have been modified to accommodate RDC restraints. Other software packages have been developed in recent years specifically for structure calculation of macromolecules from RDC data[16,24,34–38]. In addition, other analyses have been introduced to perform conformational sampling of dynamical proteins from RDC data[10,38–46]. Here we present a combination of REDCRAFT[12,16,47] and new analysis of its results to explore and assemble conformational states of a dynamical protein. A better understanding of molecular structures (including challenging proteins) and dynamics will enhance our understanding of molecular events at relevant time-scales.

## 3  Methods

The presented methodology proceeds in four conceptual steps of: structure determination, identification of onset of dynamics, classification of dynamics, and reconstruction of different conformational states. Testing and validation reported in this work is based on simulated instances of dynamics and their corresponding RDC data. We have utilized REDCRAFT and a few other auxiliary programs to achieve our objectives. The following sections detail our methodology and approach in treatment of dynamics.

### 3.1  Residual Dipolar Couplings

Residual Dipolar Couplings (RDCs) have been observed as early as 1963 [48] and have been acquired for a number of structure determination studies including small molecules [49,50], carbohydrates [19,51–53], nucleic acids [18,54–57] and proteins [3,24,58–62]. The RDC interaction phenomenon has been extensively reviewed in the literature[cite reviews] and we therefore refrain from providing an overall view of this phenomenon. The physical principles [48,63] that lead to manifestation of RDCs, and methods inducing alignment of biological macromolecules, have been fully described previously [20,64–66]. Here we briefly review those components utilized by REDCRAFT. In addition, we limit our discussion to nuclei with spin quantum number of ½, and refer to the formula in Equation 1 from which all mathematical derivation of the RDC interactions (for a pair of spin ½ nuclei) begin. In this equation, $\mu_0$ is the magnetic permeability of free space, $\gamma_i$ and $\gamma_j$ are gyromagnetic ratios of the interacting nuclei, $h$ is Planck's constant, $r$ is the distance separating nuclei $i$ and $j$, and $\theta$ is the angle between the magnetic field of the NMR device and a vector connecting atoms $i$ and $j$.

$$D_{ij} = \frac{-\mu_0 \gamma_i \gamma_j h}{(2\pi r)^3} \left\langle \frac{3\cos^2\theta_{ij}(t)-1}{2} \right\rangle \qquad (1)$$

It is important to note that the RDC value $D_{ij}$ (reported in units of Hz) is a function of the time-dependent angle $\theta(t)$ averaged over time $t$, as represented by the angular brackets in Equation 1. This time averaging phenomenon may account for molecular motions caused by: natural bond vibrations, internal dynamics, or overall tumbling of the molecule in the solution state. Mathematical transformation[48] of Equation 1 can produce a computationally amiable formulation of the RDC phenomenon, as shown in Equation 2. In this representation of the RDC interaction, $\bar{v}$ signifies the normalized orientation of the interacting vector, $s_{ij}$ denotes the $ij^{th}$ element of the Saupe order tensor matrix, $S_{ii}$ represents the principle order parameters, and $\xi$ symbolizes the Eulerian rotation matrix that relates the molecular frame to the principal alignment frame. The remaining constants have been subsumed into a single constant, $D_{max}$.

$$D = D_{max} \bar{v}^T \cdot \begin{pmatrix} s_{xx} & s_{xy} & s_{xz} \\ s_{yx} & s_{yy} & s_{yz} \\ s_{zx} & s_{zy} & s_{zz} \end{pmatrix} \cdot \bar{v} = D_{max} \bar{v}^T \cdot \xi(\alpha, \beta, \gamma) \cdot \begin{pmatrix} S_{xx} & 0 & 0 \\ 0 & S_{yy} & 0 \\ 0 & 0 & S_{zz} \end{pmatrix} \cdot \xi^T(\alpha, \beta, \gamma) \cdot \bar{v} \qquad (2)$$

Within recent years, various methods have been proposed [27,35,67–75] to estimate of the optimal order tensor. These diverse approaches exhibit some advantages over other existing methods, such as: highly accurate estimation of order tensors, estimation of order tensor or order parameters in the absence of a structure, relative order tensor estimation from unassigned RDCs, and reconstruction of interacting vectors in space from unassigned RDCs. However, the most reliable method of obtaining an order tensor is from using an assigned set of RDCs to a high-resolution structure. Due to the nature of the work presented in this paper, the estimation of order tensor will be based on assigned RDC from a structure[67,70,73] that is computed concurrent to treatment of dynamics.

### 3.2 Classification of Modes of dynamics

An informed approach to study of internal dynamics in macromolecules requires classification of different types of dynamics. This is necessary in order to select the most optimal approach to study and reconstruction of the trajectory of dynamics. To better facilitate this selection and further discussion of dynamics we enumerate three distinct dimensions of dynamics, namely: Temporal, Structural and Alignment as shown in Table 1. Along the Structural mode of dynamics we define two categories of Rigid-body and Uncorrelated modes. Similar to previous definitions[32,76–78], Rigid-body dynamics is defined as dynamical regions that maintain a constant internal structure as a function of time, while the Uncorrelated dynamics is defined as alteration of structure as a function of time. Therefore it is meaningful to describe the structure of a dynamical region if it is engaged in a Rigid-body dynamics, and not so for an Uncorrelated mode of dynamics. The temporal dimension of dynamics can be defined by two categories of Discrete-state and Continuous-state dynamics. The distinction between the two is solely based on the temporal occupancy of conformational states that are visited during the trajectory of the dynamics. The Alignment mode of dynamics can be described by homogeneous and heterogeneous modes of alignment where homogeneous mode of alignment assumes fixed alignment of the protein (within the same alignment medium) as a function of conformational changes. The heterogeneous mode permits modeling of dynamics where alignment of a protein is perturbed as a function of conformational state of the protein. Where as in the homogeneous mode of dynamics, alignment of a protein is altered as a function of the conformational changes. In principle all eight combined modes of dynamics should be possible with examples of all four combination of Structural and Temporal modes of dynamics having already been presented in the literature[6–9,79]. In this report we investigate the combination of Rigid-body, Discrete-state dynamics with the explanation that it represents biologically most likely event, and that the remaining three modes (combinations of Structural and Temporal modes) can be approximated as a Rigid-body and Discrete-state dynamics in some favorable instances.

*Table 1: Different modes of dynamics*

| Structural | Temporal | Alignment |
|---|---|---|
| Rigid-body | Continuous-state | Homogeneous |
| Uncorrelated | Discrete state-state | Heterogeneous |

The foundation of the presented work is based on reconstructing the trajectory of dynamics from discrepancies of order tensors reported from the static and dynamic domains of a protein. Therefore the first step in the study of dynamics is formulation of affects of dynamics on order tensors. Equation 3 formulates changes in the observable order tensor as a function of time (or dynamics). In this equation the variable $j$ denotes the $j^{th}$ alignment medium and integration is performed over the entire life of the dynamics. It can be argued that since biological systems perform cyclical motions (returning to some original state), therefore the lifetime of a dynamic event can be treated as finite and periodical. Discrete approximation of the continuous function shown in Equation 3 can be developed as shown in Equation 4. In this formulation $\delta t$ serves as the discrete time interval of observation, which if selected appropriately can provide an accurate approximation of a temporally continuous motion. This equation can be further simplified based on relative occupancies in different states of the dynamics. Equation 4 can be further simplified if the temporal occupancy of the conformational continuum is primarily in a small number of stable states and the contribution of the transient states is negligible. Under these conditions Equation 5 can be formulated and adopted in recovery of the primary conformational states of a Discrete-state dynamics. In this equation the entity $S^i_j$ denotes the order tensor reported from the $i^{th}$ conformational state within the $j^{th}$ alignment medium where $\rho_i$ is the relative occupancy of the $i^{th}$ state. The second constraint shown in this equation enforces the fact that the sum of all relative occupancies should equate to 1 (or 100%).

$$\widehat{S}_j = \int_{t=0}^{\infty} S_j(t)\,dt \tag{3}$$

$$\widehat{S}_j = \sum_{k=1}^{n} S_j(k \cdot \delta t) \tag{4}$$

$$\begin{cases} \widehat{S}_j = \sum_{i=1}^{n} \rho_i S^i_j \\ \text{Subject to: } \sum_{i=1}^{n} \rho_i = 1 \end{cases} \tag{5}$$

### 3.3 REDCRAFT

REDCRAFT[12,16,23,25] is designed for structure determination purely from orientational restraints. REDCRAFT deploys a powerful search mechanism that is significantly different from traditional optimization techniques, allowing for the same accuracy in recovery of structures compared to other algorithms while utilizing less data. The REDCRAFT algorithm has been previously described in depth[16,47,80]. In this section we will present only the features that are relevant during study of structure and dynamics of proteins.

REDCRAFT is well suited for the study of structure and dynamics because of its key feature of calculating the optimal structure by appending one residue at a time. This elongation process is consistent with the biological synthesis of proteins and allows for progressive examination of the rigidity assumption of a protein's structure. The averaging of order tensors due to internal dynamics results in differences in the observed order tensors between the static and rigid components of a molecule. The differences of the order tensors result in an inherent inability to produce a structure that will consistently satisfy the orientational constraints between the static and dynamical regions. These inconsistencies can be identified from REDCRAFT's *dynamic-profile*[12,47,80] that is produced during a structure calculation session. *Dynamic-profile* is formally presented and further discussed in section 3.4 and results are shown in section 4.1.

The second feature of REDCRAFT that further enables study of structure and dynamics of proteins, emanates from its ability to conduct fragmented reconstruction of a protein. In general, structure of a given protein can be created in numerous fragments because of data availability, biological importance, or study of dynamical regions that undergo Rigid-body dynamics. Study of *dynamic-profile* allows for identification of hinge regions, which can then be used to establish different dynamical domains of a protein for fragmented calculation of structures. The *dynamic-profile* has been described previously[81–83] but to facilitate better discussion, it is briefly discussed in the next section. Relevant to the current discussion, fragmented structure calculation that can be initiated from analysis of a *dynamic-profile* calculation allows structure reconstruction of all rigid components of a protein, although they may be dynamical with respect to each other. Once individual rigid components of the protein are reconstructed, they can be assembled under a dynamics scheme that will describe differences in the observed order tensors across all alignment media.

### 3.4 Dynamic Profile of REDCRAFT

The first step in investigation of dynamics of a protein is to identify the hinge regions that give rise to internal movement. It is also important to establish the structural mode of dynamics (rigid-body or uncorrelated) following discovery of the onset of dynamics. Previously presented[47,80,84] *dynamic-profile* that is produced by REDCRAFT can assist in discovery of onset of dynamics and structural mode of dynamics. An example of a typical *dynamic-profile* for a static protein is shown in Figure 1. Under typical and non-anomalous conditions, a *dynamic-profile* will start with a very low RDC-rmsd score (due to initial lack of RDC data), will monotonically increase until arriving at a maximum value that culminates around the anticipated data acquisition error, followed by a final phase that is characterized by a plateauing of the RDC-rsmd score. Any significant departure from this typical profile is indicative of some anomalous conditions. The anomalous conditions may consist of non-standard amino acid geometries (e.g. cis-Pro, impermissible dihedrals, non-standard bond lengths, etc.), existence of internal dynamics or mis-assignment of the RDCs, to name a few. Of particular interest to the discussion presented in this work, we will present alternations of *dynamic-profile* as the means to identify the onset of dynamics, and distinguish different structural modes of dynamics. Dynamic profiles can be generated for forward (N-terminus to C-terminus) or backward (C-terminus to N-terminus) analysis of the a protein. The forward and backward *dynamic-profiles* can help to corroborate the same anomalous regions with different degrees of certainty.

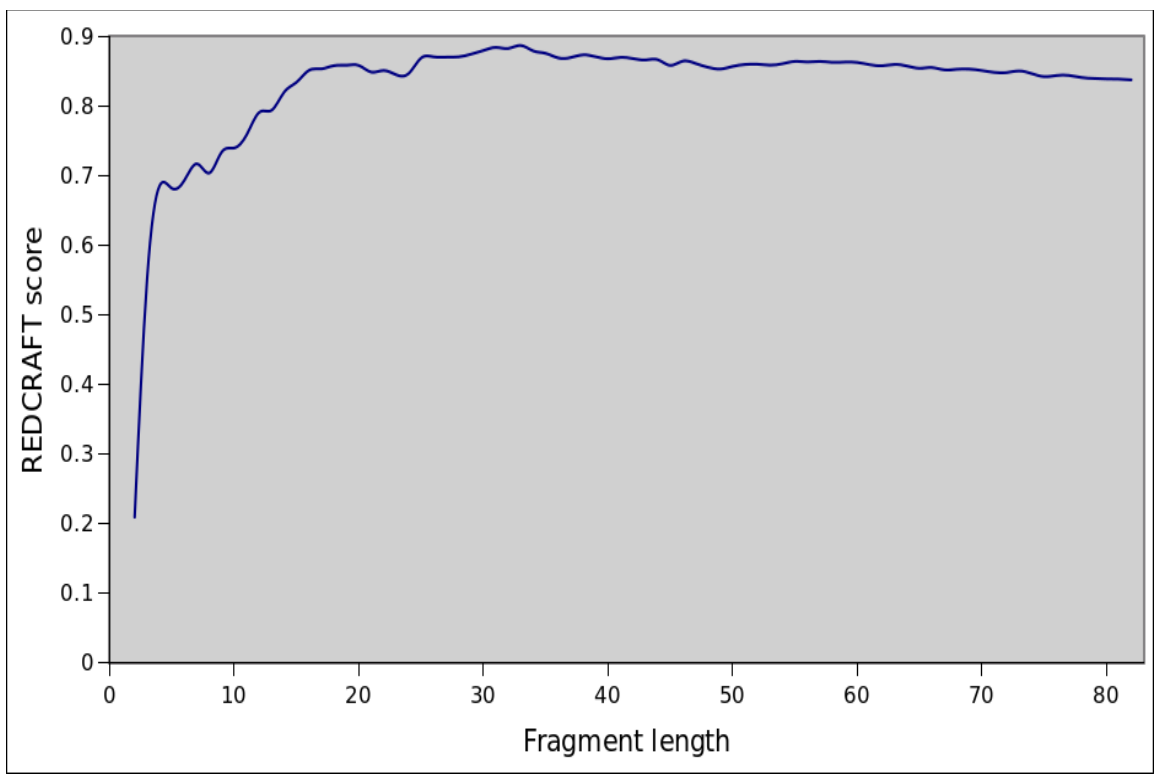

*Figure 1. Example of a typical dynamic-profile for the protein 1A1Z in the absence of internal dynamics with simulated ±1Hz of uniformly distributed noise.*

Analysis of REDCRAFT's *dynamic-profile* is based on two principles. The first principle is to identify any form of structural anomalies by observing any deviation from a typical profile. The second principle is based on the ability of REDRAFT to perform a fragmented structure determination of a protein. Once the point of anomaly is established, a new session of structure determination can be initiated a few residues in advance of the point of anomaly. The behavior of the *dynamic-profile* will be indicative of the structural mode of the dynamics. In section 4.1 we present results that demonstrate the use of this approach in discovery of onset and identification of structural mode of dynamics. Our exploration will consist of simulated Rigid-body dynamics and Uncorrelated dynamics using the protein 1A1Z as the target of our investigations. The specifics of the simulated dynamics are discussed in section 3.6.2.

### 3.5 Theoretical treatment of dynamics

The following steps describe the overall strategy to calculation of structure and characterization of dynamics that is presented in this report:

1. Proceed in structure calculation with REDCRAFT under the assumption of structural rigidity.

2. Upon identification of internal dynamics from *dynamic-profile*, embark on fragmented study of dynamics for each region that exhibits internal structural rigidity.

3. After successful completion of fragmented structure calculation, establish the rigid and dynamical fragments through comparison of observed order tensors in all alignment media. Fragments can be collected into relative rigid domains based on the similarity of their order tensors. Comparison of order tensors across different domains can establish static domains and dynamic domains.

4. Construct models of dynamics that will successfully explain the differences of the observed order tensors in all alignment media between the static and dynamic domains.

The scientific basis, technical requirements and procedures to establish steps 1-3 have been previously described and can easily be accomplished through the use of REDCRAFT and REDCAT[67,73] software packages. However, additional theoretical formulations and procedural analyses are required for step 4. To facilitate development of procedures to accomplish the objectives of step 4, we first submit that in case of a two domain dynamics it is possible to designate one of the domains as the static domain and the other as the dynamic domain. Although at first sight the principle of relative motion may appear to introduce some ambiguity to this designation, the presence of a third entity (the external magnetic field) against which all tumbling, vibrational motions and internal dynamics are observed, disambiguates the designation. It is therefore possible to uniquely designate one domain as the dynamical and the other as the static domain by simply observing the General Degree of Order[45] (GDO) of each domain. Furthermore, Equation 5 can be used as the basis of expansion to accommodate reconstruction of individual discrete states as shown in Equation 6. In this equation the term $S^a_j$ denotes the anchor order tensor in alignment medium $j$ and it signifies the order tensor that would have been observed if the dynamical domain was fixed and void of dynamics. The anchor-order tensor can be obtained from the static domain of the protein (domain with the highest GDO). The term $\xi_i$ represents the Eulerian transformation (with its three corresponding angular arguments) that maps the Rigid-body structure of the dynamical domain from any arbitrary molecular frame to the frame that defines the $i^{th}$ state of dynamics. The average observable order tensor on the left-hand side of the equation can be obtained within REDCAT by using the structure of the dynamical domain (obtained from fragmented study in REDCRAFT) and the experimentally acquired RDCs (that encompass and report the internal dynamics). Equation 6 can be used to formulate the objective function shown in Equation 7, which can be used to obtain solutions for four unknowns (relative occupancy and three Euler angles) that define each state of a given discrete dynamics. In this equation the symbol $\|.\|$ denotes magnitude of the difference-matrix by summing the square of its elements. This equation can be repeated for each alignment medium, which will contribute five additional independent equations to the overall system of equations. In total, defining *n* discrete dynamical states will require *4n-1* (relative occupancy of the last state can always be computed by one minus the sum of all other occupancies) degrees of freedom, while *m* alignment media will provide *5M* number of equations. Therefore viable solutions can be obtained so long as the criterion shown in Equation 8 is satisfied. Note that an important fact in combining information across all alignment media is that relative occupancies and orientation of the dynamical domains with respect to the static domain remain unchanged across all alignment media. We have used least-square minimization[85,86] minimization routine available in Maple 14 software package as the approach to obtain solutions to the Equation 7.

$$\begin{cases} \hat{S}_j = \sum_{i=1}^{n} \rho_i S_j^i = \sum_{i=1}^{n} \rho_i \cdot \xi(\alpha_i, \beta_i, \gamma_i) \cdot S_j^a \cdot \xi'(\alpha_i, \beta_i, \gamma_i) \\ \text{Subject to: } \sum_{i=1}^{n} \rho_i = 1 \end{cases} \quad (6)$$

$$f(\rho_{1..n}, \alpha_{1..n}, \beta_{1..n}, \gamma_{1..n}) = \begin{cases} \sum_{j=1}^{m} \left\| \hat{S} - \sum_{i=1}^{n} \rho_i \cdot \xi(\alpha_i, \beta_i, \gamma_i) \cdot S_j^a \cdot \xi'(\alpha_i, \beta_i, \gamma_i) \right\| \\ \text{Subject to: } \sum_{i=1}^{n} \rho_i = 1 \end{cases} \quad (7)$$

$$5m \geq 4n - 1 \quad (8)$$

### 3.6 Testing and Validation

Our general testing and validation strategy have relied on the use of simulated RDC data. The overall process consists of generating average sets of RDC data from different models of dynamics, reconstructing fragmented structures based on steps 1-3 as listed in section 3.5, followed by reconstructing the states from the recovered Euler rotations (after solving Equation 7). Following reconstruction of the discrete states, validation is based on quantifying the backbone deviation between the reconstructed and target states. In our experiments we utilized synthetic data from an 83 residue FADD protein (PDB ID 1A1Z) and the 213 residue human DGCR8 core (PDB ID 2YT4). The use of simulated data allows exploration in limitations of the presented methodology on theoretical models of dynamics with accessibility to the ground-truth answer, and exploration of different conditions of dynamics. The models of dynamics, summary of synthetic data and structure validation procedure used during our experiments are further described in the subsections that follow.

### 3.6.1 Simulated data

Simulation of RDC values for an arbitrary pair of nuclei requires a-priori knowledge of an order tensor. A Saupe order tensor[48] can be expressed via a 3×3 matrix, or by providing principal order parameter values $S_{xx}$, $S_{yy}$, and $S_{zz}$ and rotational Euler angles α, β, and γ. In this report we use the latter formulation of an order tensor. Using the atomic coordinates, order parameters and Euler angles; REDCAT[67,73] was used to produce computed RDC values. We have utilized a number of order tensors in our investigations as the means to passively observe the dependency of our method on order tensors. Tables 2-4 summarize the order tensors used for each of our models of dynamics.

**Table 2** Order Tensor parameters used for the simulated 2-state arc motion and the simulated DGCR8 dynamics.

|    | $S_{xx}$ | $S_{yy}$ | $S_{zz}$ | α | β | g |
|----|----------|----------|----------|---|---|---|
| S1 | 3.00×10⁻⁴ | 5.00×10⁻⁴ | -8.00×10⁻⁴ | 0º | 0º | 0º |
| S2 | -4.00×10⁻⁴ | -6.00×10⁻⁴ | 1.00×10⁻³ | 40º | 50º | -60º |

**Table 3** Order Tensor parameters used for the complex 2-state model of dynamics.

|    | $S_{xx}$ | $S_{yy}$ | $S_{zz}$ | α | β | g |
|----|----------|----------|----------|---|---|---|
| S1 | -3.00×10⁻⁴ | -5.00×10⁻⁴ | 8.00×10⁻⁴ | 0º | 0º | 0º |
| S2 | 2.00×10⁻⁴ | 5.00×10⁻⁴ | -7.00×10⁻⁴ | -40º | -50º | 60º |

**Table 4** Order Tensor parameters used for the complex 3-state model of dynamics.

|    | $S_{xx}$ | $S_{yy}$ | $S_{zz}$ | α | β | g |
|----|----------|----------|----------|---|---|---|
| S1 | 3.00×10⁻⁴ | 5.00×10⁻⁴ | -8.00×10⁻⁴ | 0º | 0º | 0º |
| S2 | 2.00×10⁻⁴ | 5.00×10⁻⁴ | -7.00×10⁻⁴ | -40º | -50º | 60º |
| S3 | -7.00×10⁻⁴ | -1.00×10⁻⁴ | 8.00×10⁻⁴ | 20º | -40º | 20º |

Simulated RDC data may also be modified to include the addition of simulated error or noise. Unless specified otherwise, all simulated RDCs are accompanied by a uniform random change in the RDC values in the range of ±1 Hz, and contain the following set of RDCs: {$C'$-$N$, $N$-$H$, $C'$-$H$, $C_\alpha$-$H_\alpha$}. To simulate different percentages of occupancies Equatoin 9 was used to average the sets of RDCs from different conformations, where $\rho_i$ and $RDC^i_n$ denote the relative occupancy and RDC values for vector $n$ in the $i^{th}$ conformational state respectively. In this equation $n$ is the total number of discrete conformational states.

$$\begin{cases} \overline{RDC}_n = \sum_{i=1}^{n} \rho_i \cdot RDC^i_n \\ \text{Subject to: } \sum_{i=1}^{n} \rho_i = 1 \end{cases} \quad (9)$$

### 3.6.2 Simulated 2-state dynamics

Our exploration of 2-state dynamics consisted of two different models of dynamics. Both models of dynamics were implemented on the FADD protein (PDB-ID 1A1Z) as an example of a helical protein. The helical nature of this protein contributed to its selection as the target of our study since helical proteins present unique challenges in study by RDC data due to the parallel orientation of their *N-H* bonds. The two models of dynamics that have been included in this report consisted of an arc motion and a more complex motion resulted from rotation about two axes. The 2-state models of arc motion were explored by rotating the $\varphi$ angle of the protein 1A1Z at the 71st residue (denoted by $\varphi_{71}$) 15°, 30° and 60°. Consequently, in the arc model of dynamics, this protein is segmented into two domains: a static domain that consists of residues 1-69 and the dynamic domain that consists of residues 73-83. An example of arc motion with 60° perturbation of $\varphi_{71}$ is shown in Figure 2. In this figure the segment of the protein illustrated in blue is the static region while the red and green domains represent the two conformations of the dynamical region. It is noteworthy that this partitioning scheme of the protein introduces additional challenges since the dynamical region is a single ideal helix.

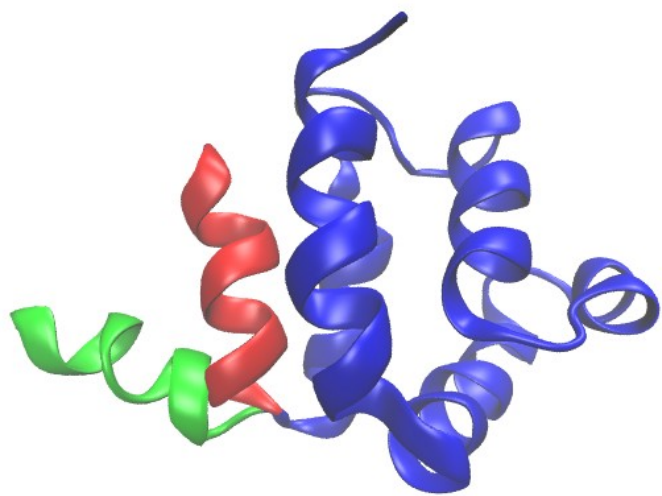

*Figure 2. 2-state arc motion of 60° at residue 71. The domain 1 of both states is seen in blue. Domain 2 of conformation 1 is seen in red and domain 2 of conformation 2 is shown in green.*

The more complex motion (example shown in Figure 3) was created by performing a 30° rotation of the $\varphi$ and $\psi$ angles at residue 58 (30° rotation of $\varphi_{58}$ followed by 30° rotation of the $\psi_{58}$) of the protein 1A1Z. In this case the two domains were defined as residues 1-56 (the static region) and residues 60-83 (dynamic region). In Figure 3 the blue portion of the structure represents the static region, while the red and orange portions of the structure represent the two alternate states of the dynamical region.

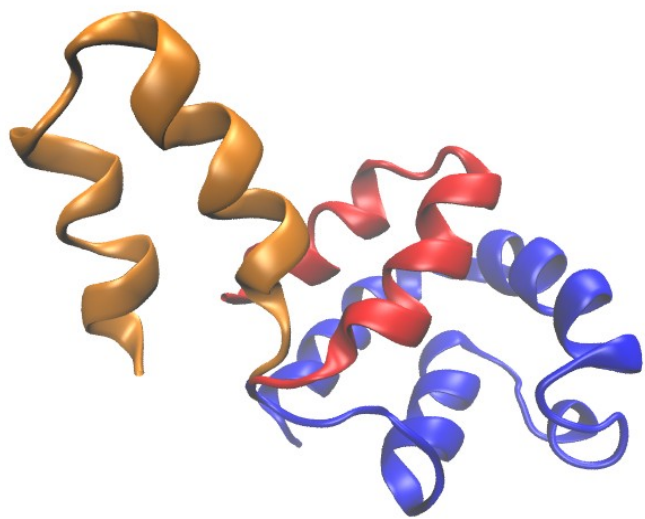

*Figure 3. 2-state complex motion created by mutating residue 58. The domain 1 of both states is seen in blue.*

*Domain 2 in conformation 1 is seen in red and domain 2 in conformation 2 is seen in orange.*

### 3.6.3 Simulated 3-state dynamics

Our exploration of 3-state dynamics consisted of building on the complex dynamics of the 2-state investigation. Here the two states from the complex 2-state motion were used as states one and two of the complex 3-state motion. The third state was created by rotating the $\varphi$ angle of residue 58 (only $\varphi_{58}$) by 60° from the original structure. As in the case of the complex 2-state motion, the domains were defined by residues 1-56 and 60-83 as the static and dynamic domains respectively. The simulated three conformations are shown in Figure 4. In this figure red, green, and orange fragments illustrate states 1, 2, and 3 of the dynamical domains while the static domain is illustrated in blue.

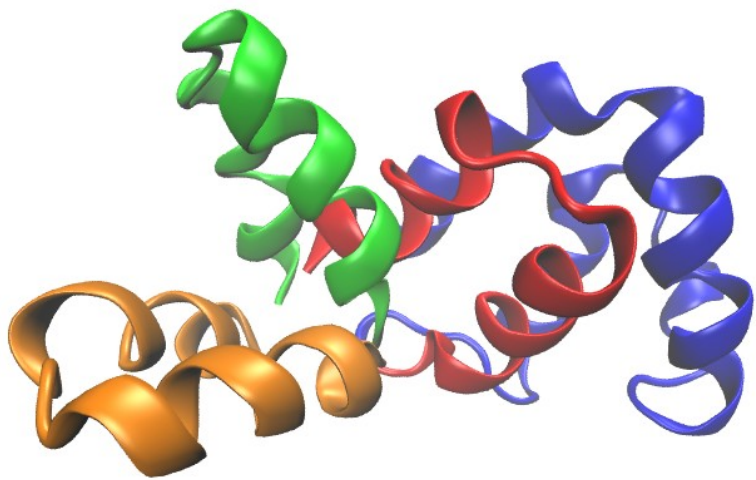

*Figure 4. Three-state complex model of dynamics. The static domain is illustrated in blue. Dynamic domain in conformational states 1, 2 and 3 are shown in red, green and orange respectively.*

### 3.6.4 Molecular dynamics simulation of DGCR8

To test our method on a more realistic motion, and one with known biological implication, we resorted to a molecular dynamic simulation of the protein DGCR8[87] to produce a plausible conformational transition. MicroRNA (miRNA) biogenesis follows a conserved succession of processing steps, beginning with the recognition and liberation of a miRNA-containing precursor (pre-miRNA) hairpin from a large primary miRNA transcript (pri-miRNA) by the Microprocessor, which consists of the nuclear RNase III Drosha and the double-stranded RNA-binding domain (dsRBD) protein DGCR8 (DiGeorge syndrome

Critical Region protein 8) which is absent in individuals with DiGeorge syndrome. DGCR8 functions as an RNA-binding protein, yet the current crystal structure of the protein would require pronounced bending of the pro-miRNA substrate to engage both distant RNA-binding domains (RBDs). Conceivably, a more likely scenario is that DGCR8 adapts[87,88] to allow for the RNA to bind by moving its two domains in tandem to create a favorable conformation to facilitate RNA binding. The motion between the two domains of DGCR8 is currently thought to be akin to a butterfly flapping its wings, with the linker region in between two RBDs as the mechanism of motion. To simulate this motion, rigid body dynamics was performed using XPLOR-NIH[32]. The published crystal structure (PDB ID 2YT4) contained several gaps (residues 497-499, 584-591, 643-648, and 702-720) which we remedied by the use of I-TASSER[89–91] structure modeling tool. The resulting structure exhibited 0.5Å of structural difference with respect to the original crystal structure with no structural gaps. Using this complete structure, 50000 steps of rigid body dynamics, with step size of 0.001psec in a 400K bath temperature were performed by keeping domain 1 (residues 17-95) and domain 2 (residues 126-203) rigid and permitting the linker region to fluctuate. The starting and ending frames of the trajectory were used as the two states with the domain 1 of both frames superimposed to create a two state jump. Using this MD simulation, DGCR8 was studied using the 2-state method. The resulting 2-state dynamics are shown as the red and green structures in Figure 5 respectively. The orientational deviation between the two dynamical states was found to be 5.4Å. In this figure, the static domain is displayed in blue and the linker region shown in purple. Average set of RDC data were computed from the two conformations using order tensors shown in Table 2 with ±0.5Hz uniformly distributed random noise added to the computed RDCs. This set of RDCs was used for reconstruction of structures by REDCRAFT in a procedure highlighted in section 3.5.

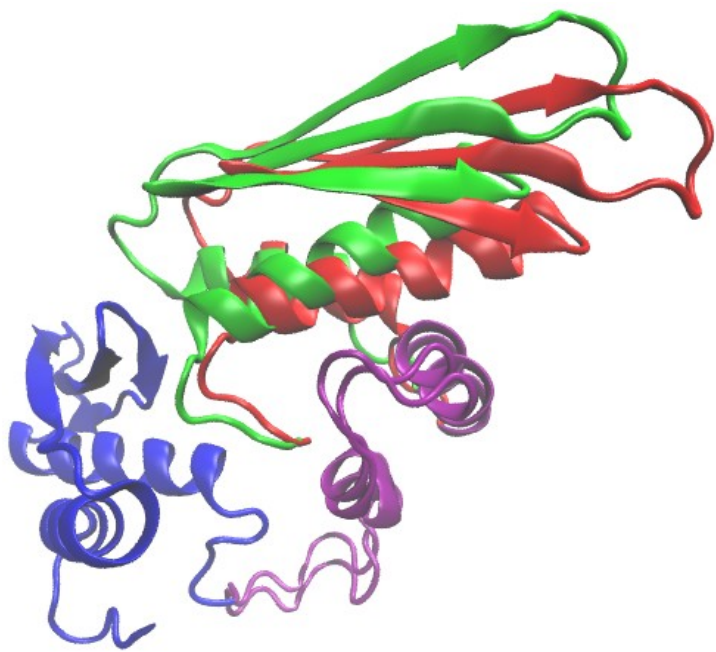

*Figure 5. Two conformations generated from Molecular Dynamics Modeling for DGCR8. Domain 1 (residues 17-95) is seen in blue, the linker region is seen in purple and the two conformations of the domain 2 are shown in red and green. The two exhibit a bb-rmsd of 5.4Å.*

### 3.6.5 Validation of results

Traditional method of reporting results for reconstructed structure of a protein is based on the measure of backbone-root-mean-square-deviation (bb-rmsd). In this report the simple use of bb-rmsd was not sufficient to report our findings since it would have generated results biased in favor of our method. Therefore a more comprehensive approach was required in order to preserve relative orientation of a protein's fragments. Complete validation of the recovered structures in this report was comprised of three consecutive steps. The first step consisted of assembling the individual structural components including the different conformations of the dynamical region. Assembly of different conformational states was accomplished by utilizing the Euler angles that were obtained from minimization of the objective function shown in Equation 7. These Euler angles facilitated the correct orientation of the conformational domains with respect to the static domain since these order tensors ($S^a_j$) were obtained from the static domain described in its molecular frame. Furthermore, because information from more than one alignment medium was used, any existing orientational degeneracies[92] (e.g. inversion degeneracy, etc.) were automatically resolved. It is important to note that upon the completion of this step, while the individual components of the protein were in correct orientational relationship with respect to each other, they may have exhibited a substantial translational error in space.

During the second step of the validation, the target structure (including all of its conformational states) was rotated to a comparable orientation with respect to the reconstructed structure to serve as a template for measurement of the bb-rmsd similarity. During this step we have used MolMol[93] visualization software to optimally superimpose the static domain of the target protein onto the static domain of the reconstructed structure through rotational and translational modifications. Completion of this step provided a measure of backbone similarity between static domains of the target and reconstructed structures. The third step of our evaluation consisted of establishing the orientational accuracy of the reconstructed conformations for the dynamic domain by allowing only translational modifications (disallowing orientational modification) of the domains. Due to the lack of existing software to perform this task we were forced to develop our own customized software for this purpose. The software "backbone" is included within the REDCRAFT software package and allows for calculation of bb-rmsd by allowing only translational modification of domains as well as optimized rotation. This allowed structural comparison of the conformational states in an orientation that was optimized based on only the static domain. It is important to note that the bb-rmsd measures that we report are an upper bound estimates and can potentially be optimized by allowing small and subtle orientational perturbation of the static domain to optimize fitness of the entire protein.

## 4 Results and Discussion

In the following sections we provide results in support of our approach in treatment of structure and dynamics of proteins. Our results first focus on the ability of REDCRAFT to accurately identify the onset of dynamics and allude to the structural mode of the dynamic. Next, we present our results in reconstruction of conformations from two and three state dynamics. We conclude our results with a discussion of limitations of the presented work and anomalies related to the study of dynamics from RDC data.

### 4.1 Discovery of onset of dynamics and structural mode of dynamics from dynamic-profile of REDCRAFT

As the first example in utility of the *dynamic-profile* we present the case of a 2-state dynamics. Here we utilized the dynamical model presented in section 3.6.2 (two states generated through perturbation of $\varphi_{71}$) and utilized the averaged RDCs to perform a forward and reverse structure calculation of the protein 1A1Z. An example of the *dynamic-profile* of a 2-state dynamic can be seen in Figure 6. In this figure the blue and red profiles correspond to the forward and reverse structure calculations respectively. As it can be seen from Figure 6, and in contrast to the typical profile shown in Figure 1, an anomalous increase has manifested in the vicinity of residue 71 on both forward and reverse sessions of REDCRAFT. This is consistent with the model of dynamics that was used during this exercise. While both forward and reverse analyses exhibit an increase in the RDC score of the *dynamic-profile*, this phenomenon is more prominently observed in the case of the reverse structure determination than the forward. This inequality arises because in the case of forward run, the anomalous region is discovered after 73 residues and RDC data from only 11 residues are inconsistent with respect to the remainder of the protein. This small portion will have a relatively smaller affect in perturbation of the RDC score reported by REDCRAFT. In contrast a much larger discrepancy is observed in the case of reverse folding of the protein because a much larger portion of the data can contribute to any observed inconsistencies.

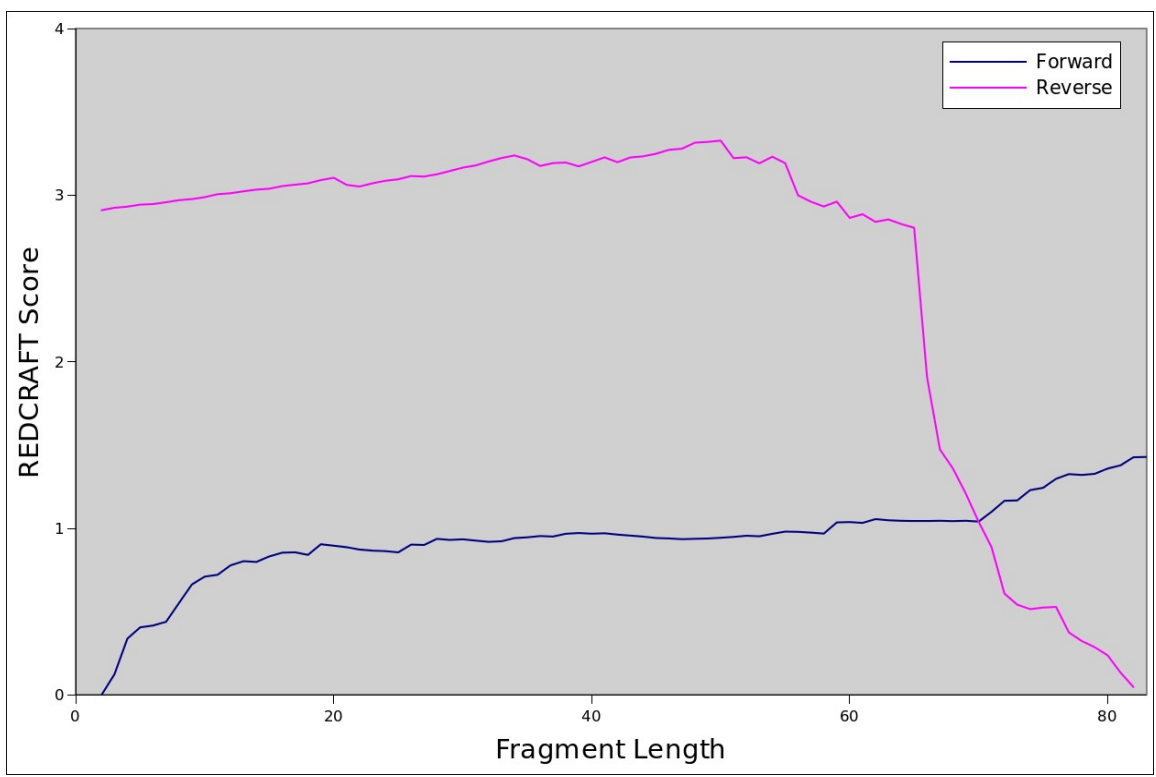

*Figure 6. An example of the dynamic profile for a 2-state model of dynamics. The blue line represents the RDC-RMSD score from REDCRAFT in forward configuration and the red line denotes REDCRAFT in reverse configuration. In this particular model of dynamics the phi angle of the 71st residue of the protein was rotated 60 degrees. The dynamic profile indicates an anomaly around that same area.*

A similar exercise was conducted for a 3-state dynamics that was described in section 3.6.3. In this model geometry of the 58th residue was altered to simulate dynamics. Figure 7 illustrates the *dynamic-profile* of this 3-state model of dynamics. Consistent with the model of dynamics, the *dynamic-profile* identifies the dynamics onset at around residue 57-58. However, unlike the previous exercise and since a larger portion of the protein is undergoing dynamics, an approximately equal increase is observed in the *dynamic-profiles* of forward and reverse structure calculation by REDCRAFT.

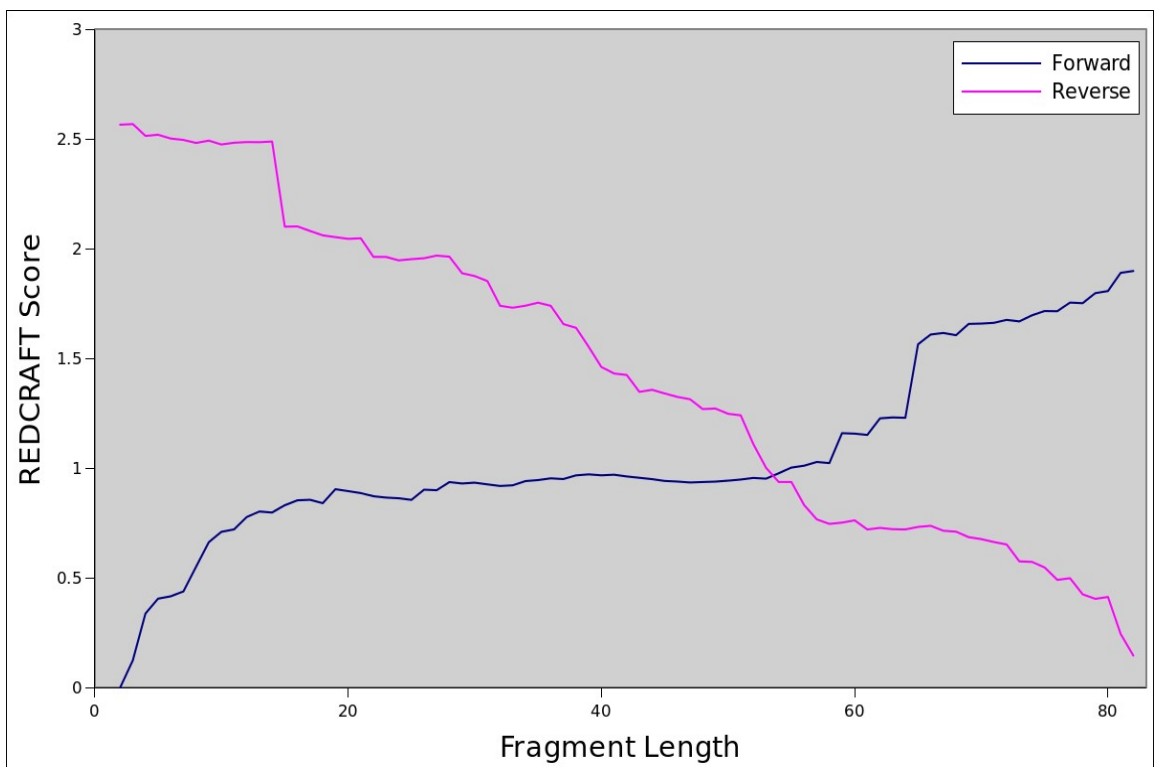

*Figure 7. An example of the dynamic profile for a 3-state model of dynamics. The blue line represents the RDC-RMSD score from REDCRAFT in forward configuration and the red line denotes REDCRAFT in reverse configuration. In this particular model of dynamics the 58th was mutated to simulate dynamics. The dynamic profile indicates an anomaly around that same area.*

The above two examples demonstrated the ability of REDCRAFT in identifying the onset of internal dynamics via the use of *dynamic-profile* analysis. Structural mode of dynamics (Rigid-body versus Uncorrelated) can be established by the use of fragmented study of a protein structure in REDCRAFT. In this context, structure calculation can be terminated prior to the onset of dynamics, and structure calculation of a new fragment can be initiated a few residues past the onset of dynamics. Analysis of the *dynamic-profile* of the new fragment can help in establishing the structural mode of dynamics. The *dynamic-profile* of the new fragment undergoing Rigid-body dynamics will exhibit a typical pattern (similar to Figure 1) since it is internally rigid and consist of a structure that is internally static as a function of time. On the other hand the Uncorrelated dynamics will exhibit a monotonically increasing score that is indicative in lack of any consistent structure as a function of time. Figures 8 and 9 illustrate examples of these two modes of dynamics simulated at $\varphi_{58}$ (in the case of rigid-body dynamics) and residues 58-83 (in the case of uncorrelated dynamics). The *dynamic-profile* of the second fragment, for the case of Rigid-body dynamics that is shown in Figure 8, exhibits a normal behavior indicating successful reconstruction of a coherent structure. D*ynamic-profile* for the case of an Uncorrelated dynamics is shown in Figure 9 and it clearly exhibits a monotonically increasing behavior that indicates absence of a coherent structure. In the case of Rigid-body dynamics, upon reconstruction of each domain's structure, a measure of relative dynamics between the two domains can be established based on comparison of their corresponding order tensors.

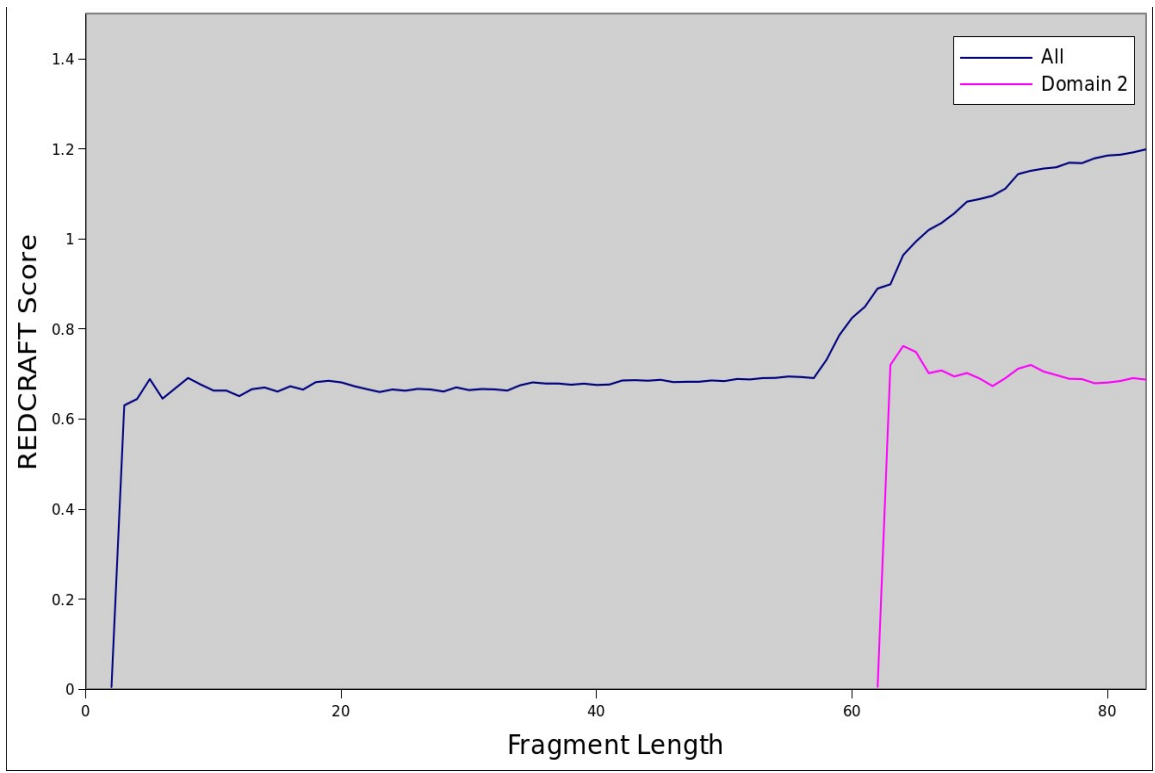

*Figure 8. Dynamic profile of discrete state dynamics.*

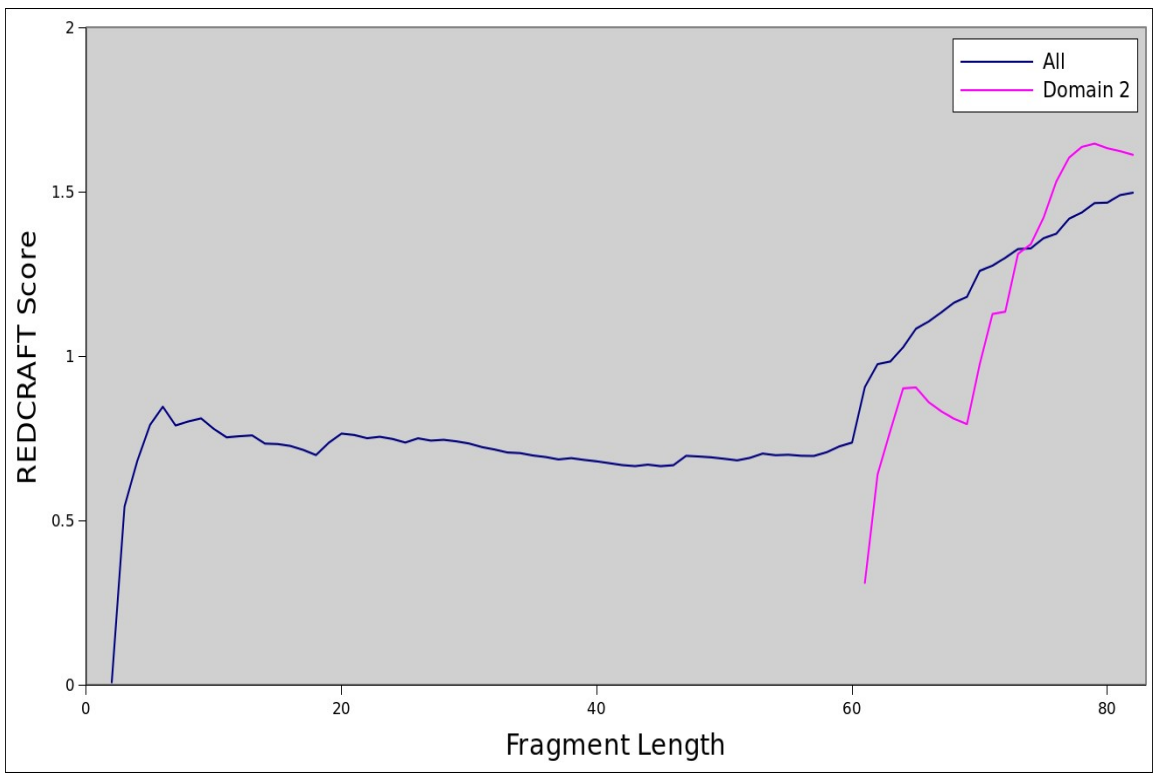

*Figure 9. Dynamic-profile of uncorrelated dynamics*

### 4.2  2-state jump dynamics

As the first step in evaluating our approach to recovery of dynamical states we resort to the arc motion of 2-state dynamics described in Section 3.6.2. Here we explored alterations of $\varphi_{71}$ by 15˚, 30˚ and 60˚ with different occupancies of the two states. We also evaluated our approach in recovery of 2-state dynamics of the complex motion (also described in section 3.6.2). Results of these experiments are shown in Tables 5, 6, 7 and 8. The results shown in these tables are segmented into sections corresponding to different states of occupancies. Our investigation of states of occupancy starts with 50% occupancy of the first state and extends into 90% in 10% increments. The row titled "Minimum" corresponds to the lowest error in the minimized objective function (shown in Equation 7) in units of Hz corresponding to an *N-H* vector. A value less than the experimental error (in this case 1Hz) indicates successful reconstruction of the domains. In this table "Conformation" denotes the two states of dynamics, "BB-RMSD" corresponds to the backbone similarity of the reconstructed states, and "Relative Occupancy" corresponds to the recovered occupancy of each state.

We begin discussion of our results with the case of 60˚ arc motion (shown in Table 5). In general all states of the dynamics were reconstructed very accurately including orientation of the states and relative occupancies. In some instances (such as 50/50) the relative occupancies were in error by as much as 13%. The only exercise that exhibited anomalous outcome was the case of 90/10. Here the first conformation was reconstructed with high degree of accuracy (0.37Å with respect to the target protein) while the second state was created with bb-rmsd of 9.4Å with respect to its corresponding state. Our explanation for this behavior is the low relative occupancy of this particular scenario marginalizes the perturbation of RDCs due to dynamics. The small perturbation of RDCs (in comparison to the noise) has therefore rendered its effect moot. This phenomenon is observed in other instances that are discussed in the future sections. Since the effect of the second state is negligible, it was reconstructed in nearly an irrelevant orientation giving rise to its high bb-rmsd to the target conformation.

*Table 5: Results for 60° arc motion.*

| | | | |
|---|---|---|---|
| 50/50 | Minimum | \multicolumn{2}{c}{$9.13\times10^{-11}$ (0.23 Hz)} |
| | Conformation | 1 | 2 |
| | BB-RMSD | 0.93Å | 1.02Å |
| | Relative Occupancy | 0.63 | 0.37 |
| 60/40 | Minimum | \multicolumn{2}{c}{$1.25\times10^{-10}$ (0.27 Hz)} |
| | Conformation | 1 | 2 |
| | BB-RMSD | 0.38Å | 0.42Å |
| | Relative Occupancy | 0.61 | 0.39 |
| 70/30 | Minimum | \multicolumn{2}{c}{$1.31\times10^{-10}$ (0.28Hz)} |
| | Conformation | 1 | 2 |
| | BB-RMSD | 0.44Å | 0.45Å |
| | Relative Occupancy | 0.72 | 0.28 |
| 80/20 | Minimum | \multicolumn{2}{c}{$2.37\times10^{-10}$ (0.37 Hz)} |
| | Conformation | 1 | 2 |
| | BB-RMSD | 0.59Å | 1.52Å |
| | Relative Occupancy | 0.85 | 0.15 |
| 90/10 | Minimum | \multicolumn{2}{c}{$2.29\times10^{-10}$ (0.37 Hz)} |

| Conformation | 1 | 2 |
|---|---|---|
| BB-RMSD | 0.37Å | 9.4Å |
| Relative Occupancy | 0.896 | 0.103 |

Next, to further investigate the sensitivity of our method with respect to the magnitude of motion we reduced the change of $\varphi_{71}$ to 30°. The results of these experiments are shown in Table 6 and are very similar to that of the 60° dynamics with the exception of the case of 80/20. In this case the second state was not able to be reconstructed with much accuracy. We suspect the reason for this inconsistency is the combination of a smaller angle of rotation and lower relative occupancy. It appears that at only 30° rotation, a state with less than 20% occupancy gets subsumed into the original state as occurred in the case of 90/10 of the previous example.

Table 6: Results for 30° arc motion

| | | | |
|---|---|---|---|
| 50/50 | Minimum | $5.12\times10^{-11}$ (0.17 Hz) | |
| | Conformation | 1 | 2 |
| | BB-RMSD | 0.46Å | 0.44Å |
| | Relative Occupancy | 0.45 | 0.55 |
| 60/40 | Minimum | $7.27\times10^{-11}$ (0.2 Hz) | |
| | Conformation | 1 | 2 |
| | BB-RMSD | 0.5Å | 0.58Å |
| | Relative Occupancy | 0.66 | 0.34 |
| 70/30 | Minimum | $1.12\times10^{-10}$ (0.26Hz) | |
| | Conformation | 1 | 2 |
| | BB-RMSD | 0.65Å | 2.00Å |
| | Relative Occupancy | 0.88 | 0.12 |
| 80/20 | Minimum | $9.4\times10^{-11}$ (0.23 Hz) | |
| | Conformation | 1 | 2 |
| | BB-RMSD | 0.66Å | 5.2Å |
| | Relative Occupancy | 0.95 | 0.05 |
| 90/10 | Minimum | $4.49\times10^{-11}$ (0.16 Hz) | |
| | Conformation | 1 | 2 |

|   | BB-RMSD | 0.5Å | 6.9Å |
|---|---|---|---|
|   | Relative Occupancy | 0.98 | 0.02 |

To further investigate the effect and behavior of our approach on small and negligible motions, we investigated the case of 15˚ arc motion. Although there exist internal variation of the second domain, the REDCRAFT dynamic profile does not identify internal dynamics indicating the absence of any anomalous behavior. Despite this finding we proceeded to reconstruct the two orientations and the results are shown in Table 7. The overarching observation that can be concluded from the results in this table consist of accurate reconstruction of the first state, and nearly complete failure to reconstruct the second state (both orientation and relative occupancy) despite achieving a low value for the objective function. This behavior is consistent with the results for 60˚ arc motion and 90/10 occupancy exercise. Both of these exercises help to establish the boundaries of the information content of the RDC data. In summary, the particular instance of 15° motion did not provide sufficient alteration of RDCs (and therefore order tensors) to indicate the existence of internal dynamics at any relative occupancies.

*Table 7: Results for 15⁰ arc motion.*

| | | | |
|---|---|---|---|
| 50/50 | Minimum | \multicolumn{2}{c}{$7.\times10^{-10}$ (0.65 Hz)} | |
| | Conformation | 1 | 2 |
| | BB-RMSD | 0.78Å | 7.8Å |
| | Rate of Occupancy | 0.96 | 0.04 |
| 60/40 | Minimum | \multicolumn{2}{c}{$1.9\times10^{-10}$ (0.34 Hz)} | |
| | Conformation | 1 | 2 |
| | BB-RMSD | 0.73Å | 9.6Å |
| | Rate of Occupancy | 0.85 | 0.15 |
| 70/30 | Minimum | \multicolumn{2}{c}{$3.48\times10^{-10}$ (0.46Hz)} | |
| | Conformation | 1 | 2 |
| | BB-RMSD | 0.68Å | 9.3Å |
| | Rate of Occupancy | 0.88 | 0.12 |

| 80/20 | Minimum | 7.13×10⁻¹⁰ (0.65 Hz) | |
|---|---|---|---|
| | Conformation | 1 | 2 |
| | BB-RMSD | 0.66Å | 9.1Å |
| | Rate of Occupancy | 0.88 | 0.12 |
| 90/10 | Minimum | 1.14×10⁻⁹ (0.82Hz) | |
| | Conformation | 1 | 2 |
| | BB-RMSD | 0.58Å | 5.2Å |
| | Rate of Occupancy | 0.95 | 0.05 |

The results from the complex 2-state model are shown in Table 8 and portray an outcome consistent with the case of arc motion. Both conformations were reconstructed with high degree of accuracy despite the complexity of the dynamics. However, it can be seen that as the relative occupancy of the second state (the ending state) decreases, the predicted orientation's BBRMSD to the target structure increases as well. This deterioration in performance is observable in the case of 80/20 and clearly so in the case of 90/10. In both cases the first state was reconstructed with reasonable accuracy while the reconstructed second state deteriorated as a function of occupancies. A relative occupancy of 10% can be seen as almost negligible in the course of a dynamic movement when using RDCs with ±1Hz of error.

*Table 8. Results for 2-state complex dynamics experiments.*

| | | | |
|---|---|---|---|
| 50/50 | Minimum | 2.27×10⁻¹⁰ (0.36 Hz) | |
| | Conformation | 1 | 2 |
| | BB-RMSD | 0.76Å | 0.83Å |
| | Rate of Occupancy | 0.42 | 0.58 |
| 60/40 | Minimum | 1.6×10⁻¹⁰ (0.31 Hz) | |
| | Conformation | 1 | 2 |
| | BB-RMSD | 1.1Å | 1.4Å |
| | Rate of Occupancy | 0.47 | 0.53 |
| 70/30 | Minimum | 1.4×10⁻¹⁰ (0.29 Hz) | |

|  | Conformation | 1 | 2 |
|---|---|---|---|
|  | BB-RMSD | 1.2Å | 1.6Å |
|  | Rate of Occupancy | 0.53 | 0.47 |
| 80/20 | Minimum | $6.04 \times 10^{-11}$ (0.19 Hz) | |
|  | Conformation | 1 | 2 |
|  | BB-RMSD | 0.69Å | 2.3Å |
|  | Rate of Occupancy | 0.66 | 0.34 |
| 90/10 | Minimum | $1.7 \times 10^{-10}$ (0.32 Hz) | |
|  | Conformation | 1 | 2 |
|  | BB-RMSD | 0.83Å | 6.33Å |
|  | Rate of Occupancy | 0.95 | 0.05 |

In summary, results in Table 7 seem to indicate that at just 15° degrees of movement, our approach is able to reconstruct one of the states (State 1) with reasonable accuracy, but fails to reconstruct the second state. However, it can be observed that when the motion is extended to a 60° or 30° movement (results shown in Tables 5 and 6 respectively), then both states can be reconstructed with reasonable accuracy so long as the relative occupancies exceed 20%. The general explanation for both cases is that the contribution of dynamics is less than the experimental noise, and therefore meaningful calculations are moot.

### 4.3  3-state jump dynamics

Results of the 3-state complex dynamics are shown in Table 9 for a number of different relative state occupancies. Similar to the case of 2-state, the relative occupancies of each exercise is listed on the left hand most column of this table. The "Minimum" value corresponds to the lowest value (in units of Hz scaled to N-H vectors) obtained from minimizing the objective function shown in Equation 7. This value helps to establish the success of the general approach; minimum values in the vicinity of the experimental noise indicate successful reconstruction of the states. The "Conformation #" denotes the conformation number, the BB-rmsd and Relative Occupancies correspond to the structural/orientational similarity and relative occupancy of each state respectively.

*Table 9. Results for 3-state dynamics experiments.*

| | | | | |
|---|---|---|---|---|
| 50/25/25 | Minimum | $2.9\times10^{-11}$ (0.13 Hz) | | |
| | Conformation # | 1 | 2 | 3 |
| | BB-RMSD | 0.95Å | 1.9Å | 0.67Å |
| | Rate of Occupancy | 0.42 | 0.32 | 0.26 |
| 34/33/33 | Minimum | $2.6\times10^{-11}$ (0.12 Hz) | | |
| | Conformation # | 1 | 2 | 3 |
| | BB-RMSD | 1.4Å | 0.38Å | 1.3Å |
| | Rate of Occupancy | 0.25 | 0.41 | 0.34 |
| 50/30/20 | Minimum | $3.4\times10^{-11}$ (0.14 Hz) | | |
| | Conformation # | 1 | 2 | 3 |
| | BB-RMSD | 1.08Å | 1.5Å | 0.4Å |
| | Rate of Occupancy | 0.32 | 0.34 | 0.34 |
| 60/30/10 | Minimum | $7.8\times10^{-11}$ (0.21 Hz) | | |
| | Conformation # | 1 | 2 | 3 |
| | BB-RMSD | 0.64Å | 1.3Å | 1.3Å |
| | Rate of Occupancy | 0.55 | 0.35 | 0.1 |

As seen in Table 9, our presented method has successfully reconstructed the conformational states and rates of occupancies with less than 2Å in structural resolution. We note that there is a higher variability in the recovered measure of relative occupancies; variations as much as 0.19%.

### 4.4 Recovery of DGCR8 discrete state dynamics

RDC data when combined with an appropriate analysis tool can provide a significant reduction in data requirements. As an example, we present results for the analysis of structure and reconstruction of a simulated 2-state model of dynamics by using only a fraction of the entire structure. In this exercise the structure calculation was limited to only a portion of the static and dynamic domains. More specifically, residues 17-41 were used as representative of the static domain, and residues 130-142 represented the dynamical region. Here we can establish the relative dynamics of two domains that are separated from each other (in space and sequence) without the need to study the entire protein. This exercise also helps to gain some insight as to the size of a fragment that is needed for successful reconstruction of dynamics. Each of the domains were reconstructed with bb-rmsd of 1.5 Å to their corresponding portion of the target protein respectively. Analysis of the order tensors strongly supported the existence of internal dynamics. Results of our conformational reconstruction are summarized in Table 10 below. Based on results shown in this table, the two states of dynamics were reconstructed with reasonable degree of accuracy. Figure 10 provides an illustration of these results for the two recovered states of DGCR8. Panels a and b of this figure depict the two conformational states of DGCR8 illustrated as transparent structures, while the opaque renderings correspond to the fragments that were reconstructed by REDCRAFT. The blue and purple regions in these pictures correspond to the static and the hinge domains respectively. The portions illustrated in red and green correspond to the two conformational states of the target protein, while the portions illustrated in yellow and orange depict the corresponding reconstructed conformations.

*Table 10. Results in recovery of DGCR8 discrete state dynamics.*

| | Minimum | $2.2 \times 10^{-10}$ (0.36 Hz) | |
|---|---|---|---|
| 50/50 | Conformation # | 1 | 2 |
| | BB-RMSD | 1.5Å | 1.5Å |
| | Rate of Occupancy | 0.54 | 0.46 |

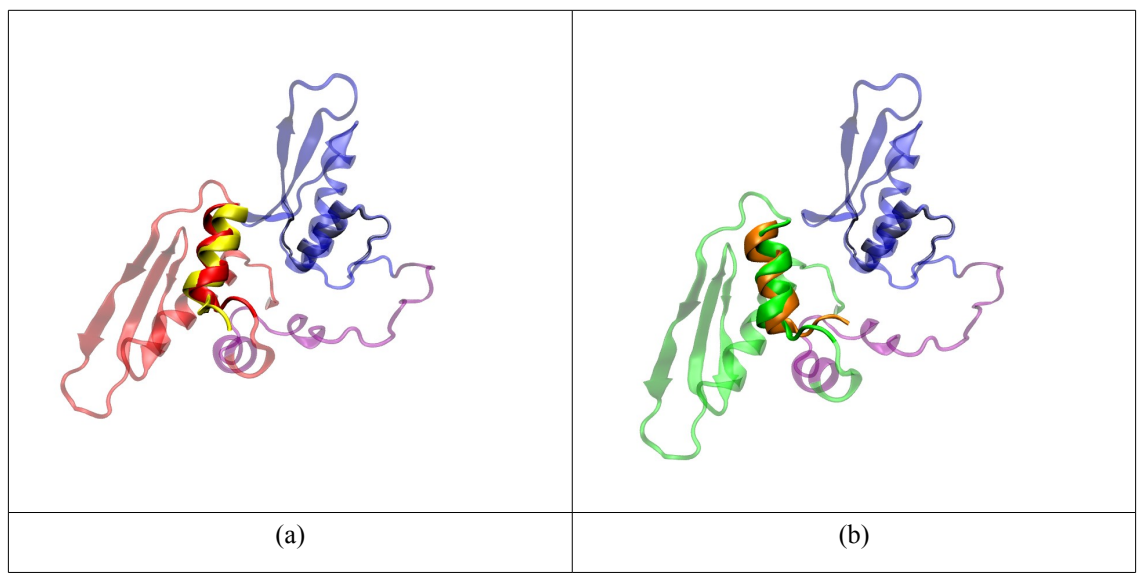

(a)        (b)

*Figure 10. Reconstructed states from the DGCR8 experiment. In both panels, the static domain is represented in blue (residues 17-41 is shown in darker blue aligned to the static domain of the target structure). On the left, the first conformation of the target structure is shown in red and the corresponding conformation as described in Table 10 is shown in yellow. On the right, the second conformation of the target protein is shown in green and the corresponding reconstructed conformation from Table 10 is shown in orange.*

### 4.5 Modeling of 2-state dynamics as 3-state dynamics or a 3-state as 2-state

Results presented in the previous section assumes a-priori knowledge of the number of dynamical states. It is reasonable to consider the cases where the number of stable conformations are not known prior to analysis. Here we present results in support of a parsimonious mechanism that assists in discovery of the appropriate jump states. More specifically, a 2-state dynamics can serve as a starting point of any investigation. Total number of conformational states can be explored incrementally until a satisfactory result is achieved. A satisfactory result is quantified by observing the fitness of the objective function defined in Equation 7 such that it equals the value of experimental error. To demonstrate this approach we resort to analysis of 3-state dynamics described in section 3.6.3. Based on this parsimonious approach, the reconstruction of conformations will proceed based on the assumption of a 2-state dynamics. Results of the 2-state recovery of the 3-state dynamics are shown in Table 11. In principle, and in agreement with the results shown in this table, the incomplete modeling should be problematic and manifest itself in an unacceptably high objective function value. Table 11 shows the results for a forced 2-state modeling of 3-state dynamics. In this table the left most column indicates the true relative occupancies of each state during the simulation of dynamics. The information marked as "Minimum" denotes fitness of the objective function (Equation 7) scaled to the units of Hz for *N-H* vectors and "Relative Occupancy" indicates the recovered values of relative occupancies. Note that in principle, and in agreement with the results shown in this table, except for the case of 60/30/10, all scenarios demonstrate a poor fitness values indicating incompleteness of a 2-state dynamics. Increasing the number of states to 3, produces the results shown in Table 9 with minimum values of the objective function that clearly indicate successful recovery of states. The case of 60/30/10 exhibited a potentially acceptable objective function because it can be treated as a two state dynamics by disregarding the state with relative occupancy of 10%. This serves as another affirmation that relative occupancies of less than 10% are negligible based on the parametrization of our simulations.

*Table 11: Results for modeling of a 3-state dynamic as a 2-state.*

| 34/33/33 | Minimum | $2.99 \times 10^{-8}$ (4.21 Hz) |

| | | | |
|---|---|---|---|
| | Relative Occupancy | 7.18×10⁻³⁰⁸ | 0.72 |
| 50/25/25 | Minimum | 4.097×10⁻⁷ (15.56 Hz) | |
| | Relative Occupancy | 0 | 1 |
| 50/30/20 | Minimum | 5.8×10⁻⁹ (1.85 Hz) | |
| | Relative Occupancy | 0.86 | 0.14 |
| 60/30/10 | Minimum | 6.9×10⁻⁹ (2.02 Hz) | |
| | Relative Occupancy | 0 | 1 |

Conversely, a 2-state dynamics can be forced to be modeled as a 3-state. In theory, a 2-state dynamics should be classified as a 3-state dynamic where two of the recovered states correspond to the two conformations, and a phantom third state with a relative occupancy of 0%. To illustrate this point we present two experiments in which we have forced 3-state recovery of a 2-state model of dynamics. Recovery of 3-state dynamics requires RDC data from at least three alignment media. In our exercise we utilized the three alignment media shown in Table 3 along with the two-state arc motion and two-state complex motion described in Section 3.6.2. In both cases equal 50% relative occupancies were used to simulate the RDC data.

*Table 12. Results for simulating 2-state dynamics in our 3-state dynamic equation.*

| | | | | |
|---|---|---|---|---|
| Arc Motion (50/50/0) | Minimum | 3.15×10⁻¹³ (0.013 Hz) | | |
| | Conformation # | 1 | 2 | 3 |
| | BB-RMSD | 0.7Å | 0.63Å | 4-7Å |
| | Rate of Occupancy | 0.47 | 0.50 | 0.03 |
| Complex Motion (50/50/0) | Minimum | 1.6×10⁻¹⁰ (0.31 Hz) | | |
| | Conformation | 1 | 2 | 3 |
| | BB-RMSD | 0.66Å | 0.6Å | 9-10Å |
| | Rate of Occupancy | 0.44 | 0.55 | 0.01 |

As can be seen from Table 12, Conformation 3 in both the arc and the complex motions have occupancy rates of 0.03 and 0.01 respectively. These conformations correspond to neither state 1 nor state 2 of their respective model of dynamics. An occupancy rate of 1-3% is, in practice, negligible making this state clearly inconsequential.

Figure 11 shows the results from the two-state arc motion with the extraneous conformation in yellow along with the two conformations (1 and 2 in Table 12) that align well with the original model of dynamics.

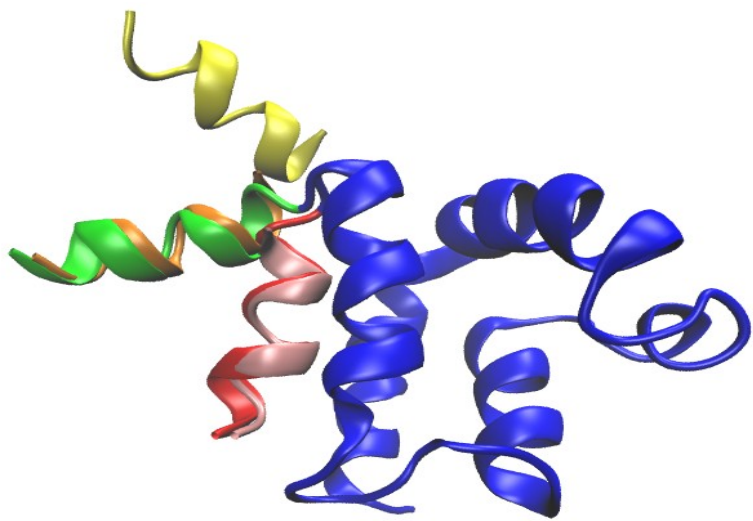

*Figure 11. The resulting structures from modeling the 2-state dynamic as a 3-state are shown here. The blue portion of the structure represents the domain 1's of each state superimposed to each other. The original domain 2 of state 1 is seen in red and domain 2 of state 2 in green. The irrelevant state is seen here in yellow, clearly not aligning to either of the two states. Conformation 1 is seen in pink (aligning to state 1) and conformation 2 is seen in orange (aligning to state 2).*

The results of these experiments are important in the sense that they reveal interesting insights into the presented method. They show that inclusion of more data improves the preciseness of the reconstruction of domains as evidenced by the low BB-RMSD of the reconstructed states in Table 12. In addition, it can be reasonably argued that in application to natural examples of dynamics where the true number of discrete states are not known, our method appears to successfully identify the proper number of states that describes a model of dynamics.

## 4.6 Limitations in recovery of discrete state dynamics

In section 4.2 we demonstrated the inability of the present work to reconstruct conformations with relative occupancies less than 10% in the cyclical life of a dynamical event. In addition, we also demonstrated the limitation in reconstructing conformational changes that are imposed by as small as 15° of arc motion. Both of these limitations are due to the overall contribution of dynamics (either relative occupancy, or small motion) relative to the experimental precision of data acquisition. Therefore in the absence of any other information, these types of limitations are universal and no approach will be able to recover useful information related to the internal dynamics.

Another category of limitations can be described that are inherent to any approach that relies on analysis of order tensors in order to recover conformational information. More specifically, these limitations arise from the fact that order tensors span a five dimensional space (degrees of freedom of an order tensor) and therefore regardless of the number of alignment media explored, no more than five independent alignment tensors can be obtained. Considering the relationship shown in Equation 8, this imposes a limitation on our approach of recovering a maximum of six conformations.

Finally, some types of dynamics are more pathological than others. There are instances of dynamics that diminish and limit the information content of RDC data and therefore impede the process of structure reconstruction. For example, previously described effects of a C3 motion[57] (three conformational states related by 120 rotation about a common axis) will alter the perceived anisotropic alignment (of any alignment medium) to appear as symmetrical with a common orientation for the $S_{zz}$ axis. This convergence of order tensors imposes some degeneracies and introduces challenges during the structure determination by RDC data. Other such pathological conditions can be demonstrated. Here we present one such example that can arise from a rotation of 90°. In this discussion, without any loss of generality we assume an arbitrary order tensor and a rotation of 90° about the z-axis with equal relative occupancies.

We demonstrate this phenomenon in the case of a 90° arc motion about $\varphi_{71}$ of the protein 1A1Z described in section 3.6.2 with relative occupancies of 50/50. Although the residual dipolar couplings were simulated with non-symmetric order tensors (in both alignment media), the observed order tensors after averaging affect of the dynamics converge to symmetrical order tensors. Figure 12 shows the SF-plots for alignment medium 1 (panel a) and medium 2 (panel b). Note that not only both order tensors have converged to a symmetrical alignment, but also the axes of symmetry have converged to the same orientation as well.

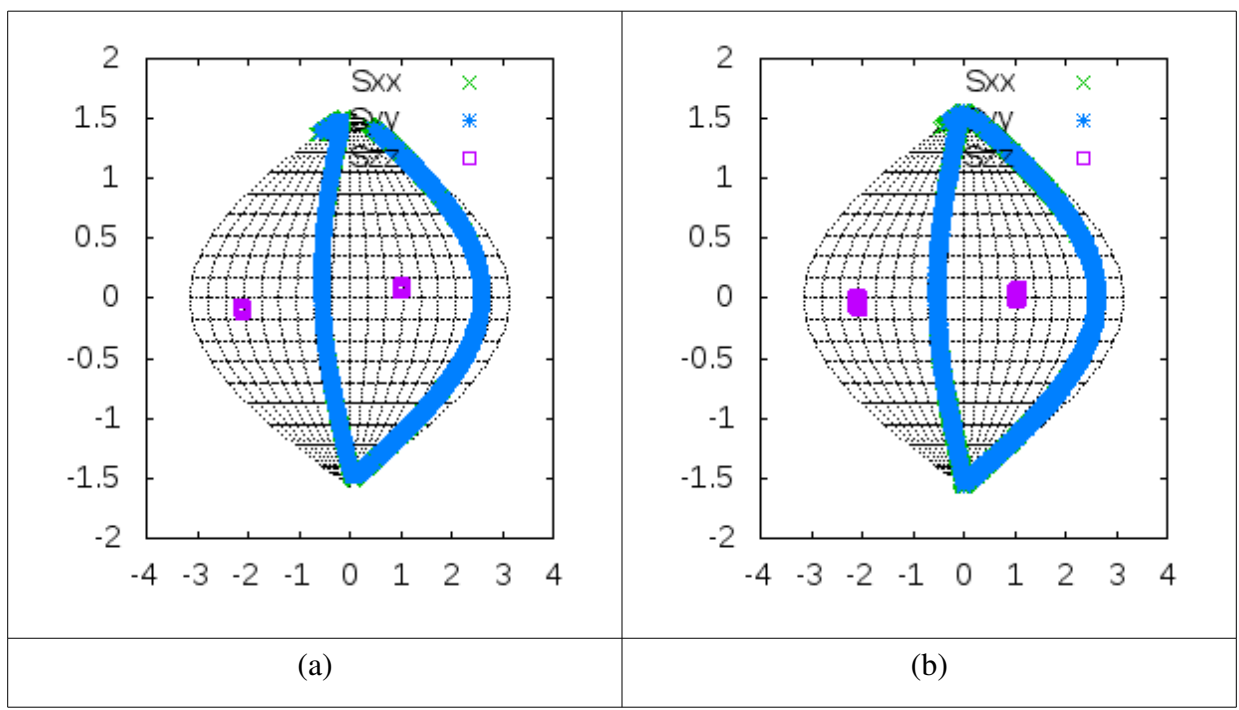

| (a) | (b) |

*Figure 12. Something here.*

## 5 Conclusions

Residual dipolar couplings can be an indispensable source of information in exploration of structure and dynamics of macromolecules. Based on the results reported here and the method that we have presented, conformational changes of internal domains can be reconstructed to within atomic resolution. In addition, the relative occupancies of each state can be reconstructed with sufficient accuracy. We have demonstrated the success of our method in application to simple dynamics that may consist of a conical rotation, rotation about two axes, and more realistic motions resulted from a Molecular Dynamic simulation. Our results demonstrated the possibility of reconstructing two and three states with high degree of accuracy. However, our investigations encountered some inherent limits in study of dynamics from RDC data. In principle, dynamics of small enough magnitude or relative occupancies, which cause perturbation of RDCs in a scale comparable or below the experimental noise can not be meaningfully reconstructed. This is consistent with principles of information theory. Based on our investigations, occupancies less than 10% of time are not sufficiently captured by RDCs for proper analysis. Our results also conclude that arc motions of less than 15º produce insufficient perturbation of RDCs for meaningful investigation of dynamics. We would however refrain from affirmatively making this particular conclusion until further studies are conducted to better explore the non-linear nature of RDC phenomenon. This non-linearity may influence the results based on selection of the orientation of the axis of rotation; rotation about one axis may cause far more alteration of RDCs than another (e.g. magic angle).

Of interest is the existence of theoretical degeneracies that could impede study of structure and dynamics with RDC data such as the example shown in section 4.6 (Figure 12). In such instances anisotropy of both alignment media converge to the same symmetrical order tensor resulting in study of structure and dynamics in effectively one alignment medium with a single symmetric order tensor. Such pathological conditions clearly interfere with study of structure and dynamics and require additional developments in order to rectify. At the current state of our investigation, such anomalies occur only with simple models of dynamics (such as arc motion) and seem to be absent in more complex cases. However, a full mathematical and theoretical investigation of this phenomenon is required. The presented formulation of dynamics can serve as a frame work for the theoretical investigation of dynamics.

We have also demonstrated the success of our method in application to as many as three conformational states. In the presence of sufficient amount of RDC data, our method can be extended for conformational analysis of as many as six states. The maximum limit of six states is imposed based on the dimensionality of order tensors. Other approaches to study of dynamics from RDCs can be envisioned that may overcome this limitation. Our future work will focus first on better understanding of degeneracies that may be encountered during the study of dynamics, and second, reformulation of the problem independent of order tensors.

Finally, we have demonstrated outcomes of our approach when the number of discrete conformational states is under or overestimated. In the case of underestimated number of conformational states, the high fitness of the objective function shown in Equation 7 will alert to the shortcomings of the assumption in number of stable states and allow for reexamination with a modified number of stable states. In the instances of overestimated number of stable states we demonstrated that the conformational states are correctly reconstructed with some additional phantom states that are represented with low relative occupancies. In this case, the phantom states can simply be discarded and the primary states utilized. The question may arise as to simply assume the maximum number of states and disregard the phantom states at the end. There are a number of reasons this strategy is discouraged. First, reconstruction of higher number of states requires RDC data from more alignment media.

Although it is always useful to acquire as much RDCs as possible, for pragmatic reasons this may not always be possible. Second, it is possible that a few phantom states (all individually negligible) may collectively represent one real state. It is therefore recommended that phantom states be eliminated one at a time followed by reevaluation with an appropriately adjusted number of states.

# 6 Acknowledgments

This work was supported by NIH Grant Numbers 1R01GM081793 and P20 RR-016461 to Dr. Homayoun Valafar.